\documentclass[thmsa,11pt]{article}
\usepackage{amsmath}
\usepackage{amssymb}

\input{tcilatex}
\begin{document}

\title{VARIATIONS ON A VISSERIAN THEME}
\author{Ali Enayat}
\maketitle

\begin{abstract}
A first order theory $T$ is \textit{tight} iff for any deductively closed
extensions $U$ and $V$ of $T$ (both of which are formulated in the language
of $T$), $U$\textit{\ }and $V$\ are bi-interpretable iff $U=V.$ By a theorem
of Visser, $\mathsf{PA}$\ (Peano Arithmetic) is tight. Here we show that $%
\mathsf{Z}_{2}$ (second order arithmetic), $\mathsf{ZF}$ (Zermelo-Fraenkel
set theory), and $\mathsf{KM}$ (Kelley-Morse theory of classes) are also
tight theories.
\end{abstract}

\section{Introduction}

The source of inspiration for this paper is located in a key result of
Albert Visser \cite[Corollaries 9.4 \& 9.6]{Albert-Tehran} concerning a
curious interpretability-theoretic feature of $\mathsf{PA}$ (Peano
arithmetic), namely:\smallskip

\noindent \textbf{1.1.~Theorem.~}(Visser) \textit{Suppose }$U$\textit{\ and }%
$V$\textit{\ are deductively closed extensions of}\textrm{\ }$\mathsf{PA}$%
\textit{\ }(\textit{both of which are formulated in the language of }$%
\mathsf{PA}$)\textit{. Then }$U$\textit{\ is a retract of }$V$\textit{\ if }$%
V\subseteq U.$\textit{\ \textbf{In particular,} }$U$\textit{\ \textbf{and} }$%
V$\textit{\ \textbf{are bi-interpretable iff} }$U=V.$\medskip

A natural reaction to Theorem 1.1 is to ask whether the exhibited
interpretability-theoretic feature of $\mathsf{PA}$ is shared by other
theories. As shown here, the answer to this question is positive\footnote{%
The $\mathsf{ZF}$-case of Theorem 1.2 was established independently in
unpublished work of Albert Visser and Harvey Friedman. I am thankful to
Albert for bringing this to my attention.}, in particular:\medskip

\noindent \textbf{1.2.~Theorem\footnote{%
See Remark 2.8 for a more complete version of this theorem.}.~}\textit{%
Theorem 1.1 remains valid if} $\mathsf{PA}$ \textit{is replaced throughout by%
} $\mathsf{Z}_{2}$ (\textit{second order arithmetic}); \textit{or by} $%
\mathsf{ZF}$ (\textit{Zermelo-Fraenkel set theory}); \textit{or by }$\mathsf{%
KM}$ (\textit{Kelley-Morse theory of classes}). \medskip 

In the remainder of this section we review some basic notions and results of
interpretability theory in order to clarify and contextualize Theorems 1.1
\& 1.2.\smallskip

\noindent \textbf{1.3.~Definitions.~}Suppose $U$ and $V$ are first order
theories, and for the sake of notational simplicity, let us assume that $U$
and $V$ are theories that \textit{support a definable pairing function }and%
\textit{\ are formulated in relational languages}. We use $\mathcal{L}_{U}$
and $\mathcal{L}_{V}$ to respectively designate the languages of $U$ and $V$%
.\medskip

\noindent \textbf{(a) }An interpretation $\mathcal{I}$ of $U$ in $V$,
written:

\begin{center}
$\mathcal{I}:U\rightarrow V$
\end{center}

\noindent is given by a translation $\tau $ of each $\mathcal{L}_{U}$%
-formula $\varphi $ into an $\mathcal{L}_{V}$-formula $\varphi ^{\tau }$
with the requirement that $V\vdash \varphi ^{\tau }$ for each $\varphi \in U$%
, where $\tau $ is determined by an $\mathcal{L}_{V}$-formula $\delta (x)$
(referred to as a \textit{domain formula}), and a mapping $P\mapsto _{\tau
}A_{P}$ that translates each $n$-ary $\mathcal{L}_{U}$-predicate $P$ into
some $n$-ary $\mathcal{L}_{V}$-formula $A_{P}$. The translation is then
lifted to the full first order language in the obvious way by making it
commute with propositional connectives, and subject to:

\begin{center}
$\left( \forall x\varphi \right) ^{\tau }=\forall x(\delta (x)\rightarrow
\varphi ^{\tau })$ and $\left( \exists x\varphi \right) ^{\tau }=\exists
x(\delta (x)\wedge \varphi ^{\tau }).$
\end{center}

\noindent Note that each interpretation $\mathcal{I}:U\rightarrow V$\/gives
rise to an \textit{inner model construction} \textit{that \textbf{uniformly}
builds a model} $\mathcal{M}^{\mathcal{I}}\models U$\/ \textit{for any} $%
\mathcal{M}\models V$. \medskip

\noindent \textbf{(b) }$U$ is \textit{interpretable} in $V$, written $%
U\trianglelefteq V$, iff there is an interpretation $\mathcal{I}%
:U\rightarrow V.$ $U$ and $V$ are \textit{mutually interpretable} when $%
U\trianglelefteq V$ and $V\trianglelefteq U.$\medskip

\noindent \textbf{(c) }We indicate the universe of each structure with the
corresponding Roman letter, e.g., the universes of structures $\mathcal{M}$, 
$\mathcal{N}$, and $\mathcal{M}^{\ast }$ are respectively $M$, $N$, and $%
M^{\ast }.$ Given an $\mathcal{L}$-structure $\mathcal{M}$ and $X\subseteq
M^{n}$ (where $n$ is a positive integer)$,$ we say that $X$ is $\mathcal{M}$-%
\textit{definable} iff $X$\ is parametrically definable in $\mathcal{M}$,
i.e., iff there is an $n$-ary formula $\varphi (x_{1},\ldots ,x_{n})$ in the
language $\mathcal{L}_{M}$ obtained by augmenting $\mathcal{L}$ with
constant symbols $\overline{m}$ for each $m\in M$ such that $X=\varphi ^{%
\mathcal{M}}$, where $\varphi ^{\mathcal{M}}=\{(a_{1},\cdot \cdot \cdot
,a_{n})\in M^{n}:\left( \mathcal{M},m\right) _{m\in M}\models \varphi (%
\overline{a}_{1},\cdot \cdot \cdot ,\overline{a}_{n})\}.$\medskip

\noindent \textbf{(d) }Suppose $\mathcal{N}$ is an $\mathcal{L}_{U}$%
-structure and $\mathcal{M}$ is an $\mathcal{L}_{V}$-structure. We say that $%
\mathcal{N}$ is \textit{parametrically} \textit{interpretable} in $\mathcal{M%
}$, written $\mathcal{N}\trianglelefteq _{\mathrm{par}}\mathcal{M}$
(equivalently: $\mathcal{M}\trianglerighteq _{\mathrm{par}}\mathcal{N}$) iff
the universe of discourse of $\mathcal{N}$, as well as all the $\mathcal{N}$%
-interpretations of $\mathcal{L}_{U}$-predicates are $\mathcal{M}$%
-definable. Note that $\trianglelefteq _{\mathrm{par}}$ is a transitive
relation.\medskip

\noindent \textbf{(e) }$U$ is a \textit{retract} of $V$ iff there are
interpretations $\mathcal{I}$ and $\mathcal{J}$ with $\mathcal{I}%
:U\rightarrow V$ and $\mathcal{J}:V\rightarrow U$, and a binary $U$-formula $%
F$ such that $F$\ is, $U$-verifiably, an isomorphism between \textrm{id}$%
_{U} $ (the identity interpretation on $U$) and $\mathcal{J}\circ \mathcal{I}
$\textsf{.} In model-theoretic terms, this translates to the requirement
that the following holds for every $\mathcal{M}\models U$:

\begin{center}
$F^{\mathcal{M}}:\mathcal{M}\overset{\cong }{\longrightarrow }\mathcal{M}%
^{\ast }:=\left( \mathcal{M}^{\mathcal{J}}\right) ^{\mathcal{I}}.$
\end{center}

\noindent \textbf{(f) }$U$ and $V$ are \textit{bi-interpretable}\footnote{%
The notion of bi-interpretability has been informally around for a long
time, but according to Hodges \cite{Hodges} it was first studied in a
general setting by Ahlbrandt and Ziegler \cite{Ahlbrandt-Ziegler}. A closely
related concept (dubbed sometimes as \textit{synonymy}, and other times as 
\textit{definitional equivalence}) was introduced by de Bouv\`{e}re \cite{de
Bouvere}. Synonymy is a stronger form of bi-interpretation; however, by a
result of Friedman and Visser \cite{Harvey-Albert}, in many cases synonymy
is implied by bi-interpretability, namely, when the two theories involved
are sequential, and the bi-interpretability between them is witnessed by a
pair of one-dimensional, identity preserving interpretations.}\textit{\ }iff
there are interpretations $\mathcal{I}$ and $\mathcal{J}$ as above that
witness that $U$ is a retract of $V$, and additionally, there is a $V$%
-formula $G,$ such that $G$\ is, $V$-verifiably, an isomorphism between $%
\mathrm{id}_{V}$ and $\mathcal{I}\circ \mathcal{J}.$ In particular, if $U$
and $V$ are bi-interpretable, then given $\mathcal{M}\models U$ and $%
\mathcal{N}\models V$, we have

\begin{center}
$F^{\mathcal{M}}:\mathcal{M}\overset{\cong }{\longrightarrow }\mathcal{M}%
^{\ast }:=\left( \mathcal{M}^{\mathcal{J}}\right) ^{\mathcal{I}}$ and $G^{%
\mathcal{N}}:\mathcal{N}\overset{\cong }{\longrightarrow }\mathcal{N}^{\ast
}:=\left( \mathcal{N}^{\mathcal{I}}\right) ^{\mathcal{J}}.$
\end{center}

\noindent We conclude this section with salient examples. In what follows $%
\mathsf{ACA}_{0}$ (arithmetical comprehension with limited induction) and $%
\mathsf{GB}$ (G\"{o}del-Bernays theory of classes) are the well-known
subsystems of $\mathsf{Z}_{2}$ and $\mathsf{KM}$ (respectively) satisfying: $%
\mathsf{ACA}_{0}$ is a conservative extension of $\mathsf{PA}$, and $\mathsf{%
GB}$ is a conservative extension of $\mathsf{ZF}$. \smallskip

\noindent \textbf{1.4.~Theorem.~}(Folklore)\textbf{\ }$\mathsf{PA}%
\trianglelefteq \mathsf{ACA}_{0}$ \textit{and }$\mathsf{ZF}\trianglelefteq 
\mathsf{GB}$; \textit{but} $\mathsf{ACA}_{0}\ntrianglelefteq \mathsf{PA}$ 
\textit{and} $\mathsf{GB}\ntrianglelefteq \mathsf{ZF}.$ \medskip

\noindent \textbf{Proof Outline.~}The first two statements have routine
proofs; the last two follow by combining (a) the finite axiomatizability of $%
\mathsf{ACA}_{0}$ and $\mathsf{GB}$, (b) the reflexivity of $\mathsf{PA}$
and $\mathsf{ZF}$ (i.e., they prove the consistency of each finite fragment
of themselves), and (c) G\"{o}del's second incompleteness theorem. \hfill $%
\square $\medskip 

\noindent By classical results of Ackermann and Mycielski, the structures $%
\left( V_{\omega },\in \right) $ and $\left( \mathbb{N},+,\cdot \right) $
are bi-interpretable, where $V_{\omega }$ is the set of hereditarily finite
sets. The two interpretations at work can be used to show Theorem 1.5 below.
In what follows, $\mathsf{ZF}_{\mathsf{fin}}$ is the theory obtained by
replacing the axiom of infinity by its negation in the usual axiomatization
of $\mathsf{ZF}$ and $\mathsf{TC}$\textit{\ }is the sentence asserting
\textquotedblleft every set has a transitive closure\textquotedblright .%
\footnote{%
More explicitly; the axioms of $\mathsf{ZF}_{\mathsf{fin}}$ consists of the
axioms of Extensionality, Empty\ Set, Pairs, Union, Power\ set, Foundation,
and$\ \lnot $Infinity, plus the scheme of Replacement\textsf{. }Note that $%
\mathsf{ZF}_{\mathsf{fin}}$ has also been used in the literature (e.g., by
the Prague school) to denote the stronger theory in which the Foundation
axiom is strengthened to the Foundation \textit{scheme}; the latter theory
is deductively identical to $\mathsf{ZF}_{\mathsf{fin}}+\mathsf{TC}$ in our
notation.} \medskip 

\noindent \textbf{1.5.~Theorem.~(}Ackermann \cite{Ackermann}, Mycielski \cite%
{Mycielski}, Kaye-Wong \cite{Kaye and Wong}) $\mathsf{PA}$\textit{\ and }$%
\mathsf{ZF}_{\mathsf{fin}}+\mathsf{TC}$ \textit{are bi-interpretable. }%
\medskip

\noindent \textbf{1.6.~Theorem.~}(E-Schmerl-Visser \cite[Theorem 5.1]{ESV}) $%
\mathsf{ZF}_{\mathsf{fin}}$ and $\mathsf{PA}$\textit{\ are not
bi-interpretable; indeed }$\mathsf{ZF}_{\mathsf{fin}}$\textit{\ is not even
a sentential retract}\footnote{%
The notion of a sentential retract is the natural weakening of the notion of
a retract in which the requirement of the existence of a definable
isomorphism between $\mathcal{M}$ and $\mathcal{M}^{\ast }$ is weakened to
the requirement that $\mathcal{M}$ and $\mathcal{M}^{\ast }$ be elementarily
equivalent.}\textit{\ of }$\mathsf{ZF}_{\mathsf{fin}}+\mathsf{TC}$\textit{.}
\bigskip

\section{Solid Theories}

The notions of solidity, neatness, and tightness encapsulated in Definition
2.1 below are only implicitly introduced in Visser's paper \cite%
{Albert-Tehran}. It is not hard to see that a solid theory is neat, and a
neat theory is tight. Hence to establish Theorems 1.1 and 1.2 it suffices to
verify the solidity of $\mathsf{PA}$, $\mathsf{Z}_{2}$, $\mathsf{ZF}$, and $%
\mathsf{KM}$. This is precisely what we will accomplish in this section. The
proof of Theorem 1.1 is presented partly as an exposition of Visser's
original proof which is rather indirect since it is couched in terms of
series of technical general lemmata, and partly because it provides a
warm-up for the proof of Theorem 2.5 which establishes the solidity of $%
\mathsf{Z}_{2}$. The proof of Theorem 2.6 establishing the solidity of $%
\mathsf{ZF}$, on the other hand, requires a brand new line of argument. The
proof of Theorem 2.7, which establishes the solidity of $\mathsf{KM}$ is the
most complex among the proofs presented here; it can be roughly described as
using a blend of ideas from the proofs of Theorems 2.5 and 2.6.\medskip

\noindent \textbf{2.1.~Definition.~}Suppose $T$ is a first order theory.

\noindent \textbf{(a)} $T$ is \textit{solid} iff the following property $%
(\ast )$ holds for all models $\mathcal{M}$, $\mathcal{M}^{\ast },$ and $%
\mathcal{N}$ of $T$:\smallskip

\noindent $(\ast )$\qquad If $\mathcal{M}\trianglerighteq _{\mathrm{par}}%
\mathcal{N}\trianglerighteq _{\mathrm{par}}\mathcal{M}^{\ast }$ and there is
an $\mathcal{M}$-definable isomorphism\textit{\ }$i_{0}:\mathcal{M}%
\rightarrow \mathcal{M}^{\ast }$, then there is an $\mathcal{M}$-definable
isomorphism\textit{\ }$i:\mathcal{M}\rightarrow \mathcal{N}$.\medskip

\noindent \textbf{(b)} $T$ is \textit{neat} iff for any two deductively
closed extensions $U$ and $V$ of $T$ (both of which are formulated in the
language of $T$), $U$\ is a retract of $V$\ if $V\subseteq U$.

\noindent \textbf{(c)} $T$ is \textit{tight} iff for any two deductively
closed extensions $U$ and $V$ of $T$ (both of which are formulated in the
language of $T$), $U$\textit{\ }and $V$\ are bi-interpretable iff $U=V.$%
\medskip

\noindent \textbf{2.1.1.~Remark.~}A routine argument shows that if $T$ and $%
T^{\prime }$ are bi-interpretable, and $T$\ is solid, then $T^{\prime }$ is
also solid. \medskip

\noindent \textbf{2.2.~Theorem.~}(Visser \cite{Albert-Tehran}) $\mathsf{PA}$ 
\textit{is solid.}\medskip

\noindent \textbf{Proof.~}Suppose $\mathcal{M}$, $\mathcal{M}^{\ast },$ and $%
\mathcal{N}$ are models of $\mathsf{PA}$ such that:

\begin{center}
$\mathcal{M}\trianglerighteq _{\mathrm{par}}\mathcal{N}\trianglerighteq _{%
\mathrm{par}}\mathcal{M}^{\ast }$, and
\end{center}

\noindent there is an $\mathcal{M}$-definable isomorphism\textit{\ }$i_{0}:%
\mathcal{M}\rightarrow \mathcal{M}^{\ast }.$ A key property\footnote{%
This important property seems to have been first noted by Feferman \cite%
{Solomon}, who used it in his proof of $\Pi _{1}$-conservativity of $\lnot 
\mathrm{Con}$($\mathsf{PA}$)\ over $\mathsf{PA.}$} of $\mathsf{PA}$ is that
if $\mathcal{M}$ is a model of $\mathsf{PA}$ and $\mathcal{N}$ is a model of
the fragment of $\mathsf{PA}$\ known as (Robinson's)\ $\mathsf{Q}$, then as
soon as $\mathcal{N}\trianglelefteq _{\mathrm{par}}\mathcal{M}$ there is an $%
\mathcal{M}$-definable initial embedding \textit{\ }$j:\mathcal{M}%
\rightarrow \mathcal{N}$, i.e., an embedding $j$ such that the image $j(%
\mathcal{M})$\ of $\mathcal{M}$ is (1) a submodel of $\mathcal{N}$, and (2)
an initial segment of $\mathcal{N}$. Hence there is an $\mathcal{M}$%
-definable initial embedding \textit{\ }$j_{0}:\mathcal{M}\rightarrow 
\mathcal{N}$ and an\textit{\ }$\mathcal{N}$-definable initial embedding 
\textit{\ }$j_{1}:\mathcal{N}\rightarrow \mathcal{M}^{\ast }$.\medskip

We claim that both $j_{0}$ and $j_{1}$ are surjective. To see this, suppose
not. Then $j(M)$ is a proper initial segment of $\mathcal{M}^{\ast }$, where 
$j$ is the $\mathcal{M}$-definable embedding \textit{\ }$j:\mathcal{M}%
\rightarrow \mathcal{M}^{\ast }$ given by $j:=j_{1}\circ j_{0}.$ But then $%
i_{0}^{-1}(j(M))$ is a proper $\mathcal{M}$-definable initial segment of $%
\mathcal{M}$ with no last element. This is a contradiction since $\mathcal{M}
$ is a model of $\mathsf{PA}$, and therefore no proper initial segment of $%
\mathcal{M}$ is $\mathcal{M}$-definable. Hence $j_{0}$ and $j_{1}$ are both
surjective; in particular $j_{0}$ serves as the desired $\mathcal{M}$%
-definable isomorphism between $\mathcal{M}$ and $\mathcal{N}$.\hfill $%
\square $\medskip

\noindent \textbf{2.2.1.~Corollary.~}$\mathsf{ZF}_{\mathsf{fin}}+\mathsf{TC}$
\textit{is solid.}\medskip

\noindent \textbf{Proof.~}In light of Remark 2.1.1, this is an immediate
consequence of coupling Theorem 2.2 and Theorem 1.5.\hfill $\square $
\medskip

Before presenting the proof of solidity of $\mathsf{Z}_{2}$ we need to state
two propositions concerning theories that prove the full scheme of induction
over some specified choice of `numbers'. The proofs of Propositions is a
straightforward adaptation of the well-known proof for the special case of $%
\mathsf{PA}$, so it is only presented in outline form.\medskip

\noindent \textbf{2.3.~Proposition.~}\textit{Let }$T$\textit{\ be a theory
formulated in a language }$\mathcal{L}$\textit{\ such that }$T$ \textit{%
interprets }$\mathsf{Q}$ \textit{via an interpretation whose domain formula
for `numbers' is} $\mathbb{N}(x)$\textit{. Furthermore assume the following
two hypotheses}:\medskip

\noindent \textbf{(a)} $T\vdash \mathsf{Ind}_{\mathbb{N}}(\mathcal{L})$, 
\textit{where} $\mathsf{Ind}_{\mathbb{N}}(\mathcal{L})$\textit{\ is the
scheme of induction over }$\mathbb{N}$ \textit{whose instances are universal
closures of }$\mathcal{L}$\textit{-formulae of the form below}:

\begin{center}
$\left( \theta (0)\wedge \forall x\left( \mathbb{N}(x)\wedge \theta
(x)\rightarrow \theta (x^{+})\right) \right) \rightarrow \forall x\left( 
\mathbb{N}(x)\rightarrow \theta (x)\right) ,$
\end{center}

\noindent \textit{where} $x^{+}$ \textit{is shorthand for the successor of} $%
x$, \textit{and} $\theta $ \textit{is allowed to have suppressed parameters;
these parameters are not required to lie in} $\mathbb{N}.$\medskip

\noindent \textbf{(b)} $\mathcal{K}\models T$\textit{\ and} $\mathcal{K}%
\trianglerighteq _{\mathrm{par}}\mathcal{N}\models \mathsf{Q}$.\textit{\
\medskip }

\noindent \textit{\textbf{Then} there is a} $\mathcal{K}$-\textit{definable
initial embedding \ }$j:\mathbb{N}^{\mathcal{K}}\rightarrow \mathcal{N}$%
.\medskip

\noindent \textbf{Proof outline.~}Since $\mathbb{N}^{\mathcal{K}}\models 
\mathsf{Q+Ind}_{\mathbb{N}}(\mathcal{L})$, the following definition by
recursion produces the desired $j$.

\begin{center}
$j(0^{\mathbb{N}^{\mathcal{K}}})=0^{\mathcal{N}}$ and $j\left( (x^{+})^{%
\mathbb{N}^{\mathcal{K}}}\right) =\left( j(x)^{+}\right) ^{\mathcal{N}}.$

\hfill $\square $
\end{center}

\noindent \textbf{2.4.~Proposition.~}\textit{Suppose} $\mathcal{K}$ \textit{%
is an }$\mathcal{L}$\textit{-structure} \textit{that interprets a model of }$%
\mathsf{Q}$ \textit{via an interpretation whose domain formula for `numbers'
is} $\mathbb{N}(x)$ \textit{and} $\mathsf{Ind}_{\mathbb{N}}(\mathcal{L})$ 
\textit{holds in} $\mathcal{K}$\textit{. Then every }$\mathcal{K}$-\textit{%
definable} \textit{proper initial segment of }$\mathbb{N}^{\mathcal{K}}$%
\textit{\ has a last element.\medskip }

\noindent \textbf{Proof.~}Easy: the veracity of $\mathsf{Ind}_{\mathbb{N}}(%
\mathcal{L})$ in $\mathcal{K}$ immediately implies that any $\mathcal{K}$%
-definable initial segment of $\mathbb{N}^{\mathcal{K}}$ with no last
element coincides with $\mathbb{N}^{\mathcal{K}}$.\hfill $\square $\medskip

\noindent \textbf{2.5.~Theorem.~}$\mathsf{Z}_{2}$ \textit{is solid.}\medskip

\noindent \textbf{Proof.} Following standard practice (as in \cite{Steve
Book}) models of $\mathsf{Z}_{2}$ are represented as two-sorted structures
of the form $\left( \mathcal{M},\mathcal{A}\right) $, where $\mathcal{M}%
\models \mathsf{PA}$, $\mathcal{A}$ is a collection of subsets of $M$, and $%
\left( \mathcal{M},\mathcal{A}\right) $ satisfies the full comprehension
scheme. Suppose $\left( \mathcal{M},\mathcal{A}\right) $, $\left( \mathcal{M}%
^{\ast },\mathcal{A}^{\ast }\right) $, and $\left( \mathcal{N},\mathcal{B}%
\right) $ are models of $\mathsf{Z}_{2}$ such that:

\begin{center}
$\left( \mathcal{M},\mathcal{A}\right) \trianglerighteq _{\mathrm{par}%
}\left( \mathcal{N},\mathcal{B}\right) \trianglerighteq _{\mathrm{par}%
}\left( \mathcal{M}^{\ast },\mathcal{A}^{\ast }\right) ,$
\end{center}

\noindent and there is an $\left( \mathcal{M},\mathcal{A}\right) $-definable
isomorphism\textit{\ }

\begin{center}
$\widehat{i_{0}}:\left( \mathcal{M},\mathcal{A}\right) \rightarrow \left( 
\mathcal{M}^{\ast },\mathcal{A}^{\ast }\right) .$
\end{center}

\noindent Note that \textit{\ }$\widehat{i_{0}}$ is naturally induced by $%
i_{0}$, where:

\begin{center}
$i_{0}:=\widehat{i_{0}}\upharpoonright M:\mathcal{M}\rightarrow \mathcal{M}%
^{\ast }$,
\end{center}

\noindent since $\widehat{i_{0}}(A)=\{i_{0}(m):m\in A\}$ for $A\in \mathcal{A%
}^{\ast }$.\medskip

\noindent It is clear that $\mathsf{Z}_{2}\vdash \mathsf{Q}^{\mathbb{N}}+%
\mathsf{Ind}_{\mathbb{N}}(\mathcal{L})$ for $\mathcal{L}:=\mathcal{L}_{%
\mathsf{Z}_{2}},$ so by Proposition 2.3, we may conclude:\medskip

\noindent (1) There is an $\left( \mathcal{M},\mathcal{A}\right) $-definable 
\textit{initial embedding }$j_{0}:\mathcal{M}\rightarrow \mathcal{N}$, and
\medskip

\noindent (2) There is an $\left( \mathcal{N},\mathcal{B}\right) $-definable 
\textit{initial embedding\ }$j_{1}:\mathcal{N}\rightarrow \mathcal{M}^{\ast
} $. \medskip

\noindent Similar to the proof of Theorem 2.2 we now argue that both $j_{0}$
and $j_{1}$ are surjective since otherwise the $\left( \mathcal{M},\mathcal{A%
}\right) $-definable embedding \textit{\ }$j:\mathcal{M}\rightarrow \mathcal{%
M}^{\ast }$ given by $j:=j_{1}\circ j_{0}$ will have the property that $j(M)$
is a proper initial segment of $\mathcal{M}^{\ast }$, which in turn implies
that $i_{0}^{-1}(j(M))$ is an $\left( \mathcal{M},\mathcal{A}\right) $%
-definable proper initial segment of $\mathcal{M}$ with no last element,
which contradicts Proposition 2.4. Hence (1) and (2) can be strengthened
to:\medskip

\noindent (1$^{+}$) There is an $\left( \mathcal{M},\mathcal{A}\right) $%
-definable \textit{isomorphism }$k_{0}:\mathcal{N}\rightarrow \mathcal{M}$,
and \medskip

\noindent (2$^{+}$) There is an $\left( \mathcal{N},\mathcal{B}\right) $%
-definable \textit{isomorphism\ }$k_{1}:\mathcal{M}^{\ast }\rightarrow 
\mathcal{N}$. \medskip

\noindent Let $\widehat{k_{0}}:\left( \mathcal{N},\mathcal{B}\right)
\rightarrow \left( \mathcal{M},\mathcal{A}\right) $ be the natural extension
of $k_{0}$, i.e., $\widehat{k_{0}}(n):=k_{0}(n$) for $n\in N$, and $\widehat{%
k_{0}}(B)=\{k_{0}(n):n\in B\}$ for $B\in \mathcal{B}$. Note that the $\left( 
\mathcal{M},\mathcal{A}\right) $-definability of $k_{0}$, along with the
veracity of the comprehension scheme in $\left( \mathcal{M},\mathcal{A}%
\right) $ assures us that $\widehat{k_{0}}(B)\in \mathcal{A}$ for each $B\in 
\mathcal{B}$. Therefore $\widehat{k_{0}}$ is an \textit{embedding.} Using an
identical reasoning, since $\left( \mathcal{M}^{\ast },\mathcal{A}^{\ast
}\right) $ is $\left( \mathcal{N},\mathcal{B}\right) $-definable by
assumption, we can extend $k_{1}$ to an embedding $\widehat{k_{1}}:\left( 
\mathcal{M}^{\ast },\mathcal{A}^{\ast }\right) \rightarrow \left( \mathcal{N}%
,\mathcal{B}\right) $. Let $\widehat{k}:=\widehat{k}_{0}\circ \widehat{k}%
_{1}\circ \widehat{i}_{0}$. Then:\medskip

\noindent (3) $\widehat{k}:\left( \mathcal{M},\mathcal{A}\right) \rightarrow
\left( \mathcal{M},\mathcal{A}\right) $ and $\widehat{k}$ is an $\left( 
\mathcal{M},\mathcal{A}\right) $-definable embedding.\medskip

\noindent The proof of Theorem 2.5 will be complete once we verify that $%
\widehat{k_{0}}$ is surjective. Since we already know that $k_{0}$ is
surjective, it suffices to check that $\mathcal{A}=\widehat{k_{0}}(\mathcal{B%
}):=\{\widehat{k_{0}}(B):B\in \mathcal{B}\}.$ Observe that the restriction $%
k:$ $\mathcal{M}\rightarrow \mathcal{M}$ of $\widehat{i}$ to `numbers' is an
automorphism of $\mathcal{M}$, thanks to (1$^{+}$), (2$^{+}$), and the
assumption that $\widehat{i_{0}}$ is an isomorphism. But since $\widehat{k}$
is $\left( \mathcal{M},\mathcal{A}\right) $-definable $i(m)=m$ for all $m\in
M$, thanks to the veracity of $\mathsf{Ind}_{\mathbb{N}}(\mathcal{L})$ in $%
\left( \mathcal{M},\mathcal{A}\right) $, for $\mathcal{L=L}_{\mathsf{Z}%
_{2}}, $ which in turn implies that $\widehat{k}$ is just the identity
automorphism on $\left( \mathcal{M},\mathcal{A}\right) $. Hence $\widehat{%
k_{0}}$ and $\widehat{k_{1}}$ are both surjective.\hfill $\square $\medskip

In the following corollary, $\widetilde{\mathsf{ZF}}$ is the result of
substituting the Replacement scheme in the usual axiomatization of $\mathsf{%
ZF}$ (e.g., as in \cite{Kunen Text}) with the scheme of Collection, whose
instances consist of universal generalizations of formulae of the form $%
\left( \forall x\in a\ \exists y\ \varphi (x,y)\right) \rightarrow \left(
\exists b\ \forall x\in a\ \exists y\in b\ \varphi (x,y)\right) $, where the
parameters of $\varphi $ are suppressed.\medskip

\noindent \textbf{2.5.1.~Corollary.~}\textit{The following theory} $T$ 
\textit{is solid}:

\begin{center}
$T:=$ $\widetilde{\mathsf{ZF}}\backslash \{\mathsf{Power}\ \mathsf{Set}\}$ + 
$\forall x\ \left\vert x\right\vert \leq \aleph _{0}.$
\end{center}

\noindent \textbf{Proof.~}In light of Remark 2.1.1, this is an immediate
consequence of Theorem 2.5 and the well-known bi-interpretability of $T$
with $\mathsf{Z}_{2}+\Pi _{\infty }^{1}$-$\mathsf{AC}$, where $\Pi _{\infty
}^{1}$-$\mathsf{AC}$ is the scheme of choice.\footnote{%
This bi-interpretability was first explicitly noted by Mostowski in the
context of the so-called $\beta $-models of $\mathsf{Z}_{2}+\Pi _{\infty
}^{1}$-$\mathsf{AC}$ (which correspond to well-founded models of $T$). See 
\cite[Theorem VII.3.34]{Steve Book} for a refined version of this
bi-interpretability result.}\hfill $\square $\medskip

\noindent \textbf{2.6.~Theorem.~}$\mathsf{ZF}$ \textit{is solid.}\medskip

\noindent \textbf{Proof.~}Suppose $\mathcal{M}$, $\mathcal{M}^{\ast },$ and $%
\mathcal{N}$ are models of $\mathsf{ZF}$ such that:

\begin{center}
$\mathcal{M}\trianglerighteq _{\mathrm{par}}\mathcal{N}\trianglerighteq _{%
\mathrm{par}}\mathcal{M}^{\ast }$,
\end{center}

\noindent and there is an $\mathcal{M}$-definable isomorphism\textit{\ }$%
i_{0}:\mathcal{M}\rightarrow \mathcal{M}^{\ast }$. Since $\mathcal{M}$
injects $M$ into $M^{\ast }$ via $i_{0}$, and $M^{\ast }\subseteq N$, we
have:$\medskip $

\noindent (1) $N$ is a proper class as viewed from $\mathcal{M}$. $\medskip $

\noindent Let $E:=\ \in ^{\mathcal{M}^{\ast }}.$ $E$ is both extensional and
well-founded as viewed from $\mathcal{N}$; extensionality trivially follows
from the assumption that $\mathcal{M}^{\ast }\models \mathsf{ZF}$, and
well-foundedness can be easily verified using the assumptions that $\mathcal{%
M}\trianglerighteq _{\mathrm{par}}\mathcal{N}$ and $i_{0}$ is an $\mathcal{M}
$-definable isomorphism between $\mathcal{M}$ and $\mathcal{M}^{\ast }$. We
wish to show that $E:=\ \in ^{\mathcal{M}^{\ast }}$ is \textit{set-like}%
\footnote{%
In the context of $\mathsf{ZF}$, the extension of a binary formula $R(x,y)$
is set-like iff for every set $s$ there is a set $t$ such that $%
t=\{x:R(x,s)\}.$}\textit{\ }as viewed from $\mathcal{N}$, i.e., for every $%
c\in M^{\ast },$ $c_{E}$ is a set (as opposed to a proper class) of $%
\mathcal{N}$, where $c_{E}:=\{x\in M^{\ast }:xEc\}.$ This will take some
effort to establish. We will present the argument in full detail, especially
because a natural adaptation of the same argument will also work in one of
the stages of the proof of Theorem 2.7 (establishing the solidity of $%
\mathsf{KM}$), and will therefore be left to the reader. We will first show
that $E$ is set-like when restricted to $\mathbf{Ord}^{\mathcal{M}^{\ast }}$.%
\footnote{%
Note that if the axiom of choice holds in $\mathcal{M}^{\ast },$ then by
Zermelo's well-ordering theorem, from the point of view of $\mathcal{N}$\
the set-likeness of $E$ when restricted to $\mathbf{Ord}^{\mathcal{M}^{\ast
}}$ immediately implies the set-likeness of $E$.} To this end, let $\delta
\in \mathbf{Ord}^{\mathcal{M}^{\ast }}$, $\delta _{E}:=\{m\in M^{\ast
}:mE\delta \}$, and consider the $\mathcal{N}$-definable ordered structure $%
\Delta :$

\begin{center}
$\Delta :=\left( \delta _{E},E\cap \delta _{E}^{2}\right) $.
\end{center}

\noindent It is clear, thanks to $i_{0}$, that $\mathcal{N}$ views $\Delta $
as a \textit{well-founded} linear order in the strong sense that every
nonempty $\mathcal{N}$-definable subclass of $\delta _{E}$ has an $E$-least
member. In particular, $\Delta $ is a linear order in which every element
other than the last element (if it exists) has an immediate successor. Given 
$\mathcal{\gamma }\in \mathbf{Ord}\mathcal{^{\mathcal{N}}}$ let $o(\Delta 
\mathcal{)\geq \gamma }$ be an abbreviation for the statement:

\begin{center}
\textquotedblleft there is some set $f$ such that $f$ is the (graph of) an
order preserving function between $(\mathcal{\gamma },\in )$ and an initial
segment of $\Delta $\textquotedblright $,$
\end{center}

\noindent and let $o(\Delta \mathcal{)}\geq \mathbf{Ord}$ abbreviate
\textquotedblleft $\forall \gamma \in \mathbf{Ord\ }o(\Delta \mathcal{)\geq
\gamma }$\textquotedblright . We wish to show that the statement $o(\Delta 
\mathcal{)}\geq \mathbf{Ord}$ does not hold in $\mathcal{N}.$ Suppose it
does. Then arguing in $\mathcal{N},$ for each $\gamma \in \mathbf{Ord}$
there is an order-preserving map $f_{\gamma }$ which embeds $\left( \gamma
,\in \right) $ onto an initial segment of $\Delta .$ Moreover, such an $%
f_{\gamma }$ is unique since it is a theorem of $\mathsf{ZF}$ that no
ordinal has a nontrivial automorphism. Hence if $\gamma \in \gamma ^{\prime
} $, then $f_{\gamma }\subseteq f_{\gamma ^{\prime }}$ and therefore $%
f:=\cup \left\{ f_{\gamma }:\gamma \in \mathbf{Ord}\right\} $ serves as an
order-preserving $\mathcal{N}$-definable injection of $\mathbf{Ord}\mathcal{%
^{\mathcal{N}}}$ onto an initial segment of $\Delta $. Invoking the
assumption $\mathcal{M}\trianglerighteq _{\mathrm{par}}\mathcal{N}$ this
shows that $\mathcal{M}$ must view $\mathcal{N}$ as well-founded because the
map $\rho ^{\mathcal{N}}:$ $(N,\in ^{\mathcal{N}})\rightarrow \left( \mathbf{%
Ord,\in }\right) ^{\mathcal{N}}$ (where $\rho $ is the usual rank function)
is $\in ^{N}$-preserving and $\mathcal{N}$-definable, and therefore $%
\mathcal{M}$-definable since $\mathcal{M}\trianglerighteq _{\mathrm{par}}%
\mathcal{N}$. This allows us to conclude that:$\medskip $

\noindent (2) $\mathcal{M}$ views $(N,\in ^{\mathcal{N}})$ as a well-founded
extensional structure of ordinal height at most $i_{0}^{-1}(\delta )\in 
\mathbf{Ord}\mathcal{^{\mathcal{M}}}$.$\medskip $

At this point we wish to invoke an appropriate form of Mostwoski's collapse
theorem in order to show that (2) implies that $N$ is a \textit{set} from
the point of view of $\mathcal{M}$. To this end, consider $\mathsf{KP}$
(Kripke-Platek set theory) whose axioms consist of Extensionality, Empty
Set, Pairs, Union, $\Pi _{1}$-Foundation, and $\Sigma _{0}$-Collection%
\footnote{%
It is well-known that $\Sigma _{1}$-Collection is provable in $\mathsf{KP}$,
which enables $\mathsf{KP}$\ to carry out\textsf{\ }$\Sigma _{1}$%
-recursions. Also note that the formulation of $\mathsf{KP}$ in many
references (including Barwise's monograph \cite{Barwise}) that focus on
admissible set theroy includes the full scheme of Foundation since
admissible sets are transtive and automatically satisfy $\Pi _{\infty }$%
-Foundation. Our forumlation of $\mathsf{KP}$ is taken from Mathias' paper 
\cite{Adrian}.}. It is well-known that $\mathsf{KP}$\ is finitely
axiomatizable, and that, provably in \textsf{KP}, $\rho $ (the rank
function) is an $\in $-homomorphism of the universe onto the class $\mathbf{%
Ord}$ of ordinals. Let $\mathsf{KPR}$ (Kripke-Platek set theory with ranks)
be the strengthening of $\mathsf{KP}$ with the axiom that states that $%
\left\{ x:\rho (x)<\alpha \right\} $ is a set for each $\alpha \in \mathbf{%
Ord}.$ Theorem 2.6.1 below can be either seen as a scheme of theorems of $%
\mathsf{ZF}$, or a single theorem of G\"{o}del-Bernays theory of classes. $%
\medskip $

\noindent \textbf{2.6.1.~Theorem.~}\textit{If} $\mathsf{KPR}$ \textit{holds
in} $\mathcal{N}$, \textit{and} $\mathbf{Ord}\mathcal{^{\mathcal{N}}}\cong
\alpha \in \mathbf{Ord}$, \textit{then }$\mathcal{N}$ \textit{is isomorphic
to a transitive substructure of} $(V_{\alpha },\in ).\medskip $

\noindent \textbf{Proof outline.} Let $h:$ $\alpha \rightarrow \mathbf{Ord}%
\mathcal{^{\mathcal{N}}}$ witness the isomorphism of $\alpha $ and $\mathbf{%
Ord}\mathcal{^{\mathcal{N}}}$, and for $\gamma <\alpha $ let $\mathcal{N}%
_{\gamma }:=(V_{h(\gamma )},\in )^{\mathcal{N}}.$ A routine induction on $%
\gamma <\alpha $ shows that there is a unique embedding $j_{\gamma }:%
\mathcal{N}_{\gamma }\rightarrow (V_{\gamma },\in )$ whose range is
transitive. This implies that if $\delta <\gamma <\alpha $, then $j_{\delta
}\subseteq j_{\gamma }.$ It is then easy to verify that $j_{\alpha }:%
\mathcal{N}\rightarrow (V_{\alpha },\in )$ is an embedding with a transitive
range, where $j_{\alpha }:=\cup \left\{ j_{\gamma }:\gamma <\alpha \right\}
. $\hfill $\square $\medskip \textit{\ }

\noindent By coupling (2) with Theorem 2.6.1 we can conclude that $N$ forms
a set in $\mathcal{M}$, thus contradicting (1). This concludes our
verification of the failure of $o(\Delta \mathcal{)}\geq \mathbf{Ord}$
within $\mathcal{N}$.\medskip

\noindent The failure of $o(\Delta \mathcal{)}\geq \mathbf{Ord}$ in $%
\mathcal{N}$ allows us to choose $\mathcal{\gamma }_{0}\in \mathbf{Ord}%
\mathcal{^{\mathcal{N}}}$ such that $\mathcal{N}$ views $\mathcal{\gamma }%
_{0}$ to be the first ordinal $\mathcal{\gamma }$ such that $o(\Delta 
\mathcal{)\geq \gamma }$ is false. We claim that $\mathcal{\gamma }_{0}$ is
a successor ordinal of $\mathbf{Ord}\mathcal{^{\mathcal{N}}}$. If not, then,
arguing in $\mathcal{N}$, for each $\beta \in \mathcal{\gamma }_{0}$ there
is a unique order-preserving map $f_{\beta }$ which maps $\left( \beta ,\in
\right) $ onto an initial segment of $\Delta ,$ and $f_{\beta }\subseteq
f_{\beta ^{\prime }}$ whenever $\beta \in \beta ^{\prime }\in \mathcal{%
\gamma }_{0}$, then $f_{\beta }\subseteq f_{\beta ^{\prime }}$. Therefore $%
\cup \left\{ f_{\gamma }:\gamma \in \mathcal{\gamma }_{0}\right\} $ serves
as an order-preserving map between $\left( \mathcal{\gamma }_{0},\in \right) 
$ and an initial segment of $\Delta $, contradicting the choice of $\mathcal{%
\gamma }_{0}.$ Hence $\mathcal{\gamma }_{0}=\beta _{0}+1$ for some $\mathcal{%
\beta }_{0}\in \mathbf{Ord}\mathcal{^{\mathcal{N}}}.$ This makes it clear
that:$\medskip $

\noindent (3) $f_{\beta _{0}}$ is a bijection between $\beta _{0}$ and $%
\delta _{E}$, $\medskip $

\noindent since if the range of $f_{\beta _{0}}$ is not all of $\delta _{E},$
then the range $\mathrm{ran}(f_{\beta _{0}})$ of $f_{\beta _{0}}$ is a
proper initial segment of $\Delta $, and $f_{\beta _{0}}$ could be extended
to an order-preserving map $f_{\mathcal{\gamma }_{0}}$ with domain $\mathcal{%
\gamma }_{0}$ by setting:

\begin{center}
$f_{\mathcal{\gamma }_{0}}(\beta _{0})=\min (\delta _{E}\backslash \mathrm{%
ran}(f_{\beta _{0}})).$
\end{center}

\noindent Thanks to (3), we now know that, as viewed by $\mathcal{N}$, $E$
is set-like when restricted to $\mathbf{Ord}^{\mathcal{M}^{\ast }}.$ To
verify the set-likeness of $E$ in $\mathcal{N}$ it is sufficient to show
that $s_{E}$ forms a set in $\mathcal{N}$, where $s_{E}:=\{m\in M^{\ast
}:mEs\}$ and $s:=V_{\delta }^{\mathcal{M}^{\ast }}$ for some $\delta \in 
\mathbf{Ord}^{\mathcal{M}^{\ast }}$ such that $\mathsf{KPR}$ holds in $%
V_{\delta }^{\mathcal{M}^{\ast }}$, since such ordinals $\delta $ are
cofinal in $\mathbf{Ord}^{\mathcal{M}^{\ast }}$ by the Reflection Theorem of 
$\mathsf{ZF}$. Consider the $\mathcal{N}$-definable structure $\Sigma :$

\begin{center}
$\Sigma :=\left( s_{E},E\cap s_{E}^{2}\right) $.
\end{center}

\noindent Since $\Sigma $ is a model of $\mathsf{KPR}$ whose set of ordinals
is isomorphic to $\beta _{0}$, by Theorem 2.6.1 (applied within $\mathcal{N}$%
) there is an $\mathcal{N}$-definable embedding of $\Sigma $ onto a
(transitive) subset of $V_{\beta _{0}}^{\mathcal{N}}$. This makes it evident
that $s_{E}$ forms a set in $\mathcal{N}$. Combined with (2) this allows us
to conclude:\medskip

\noindent (4) $E$ is extensional, set-like, and well-founded within $%
\mathcal{N}$. \medskip

\noindent At this point we invoke the Class-form of Mostowski's Collapse
Theorem:\textbf{\medskip }

\noindent \textbf{2.6.2.~Theorem.~}\cite[Theorem 5.14]{Kunen Text} \textit{%
Suppose }$E$\textit{\ is a well-founded, set-like class, and extensional on
a class }$M^{\ast }$\textit{; then there is a transitive class }$S$\textit{\
and a 1-1 map }$G$\textit{\ from }$M^{\ast }$\textit{\ onto }$S$\textit{\
such that }$G$\textit{\ is an isomorphism between }$\left( M^{\ast
},E\right) $\textit{\ and }$\left( S,\mathbf{\in }\right) $\textit{.}\textbf{%
\medskip }

\noindent Theorem 2.6.2 together with (4) assure us of the existence of an $%
\mathcal{N}$-definable $S\subseteq N$ such that $S$ is transitive from the
point of view of $\mathcal{N}$, and which has the property that there is an $%
\mathcal{N}$-definable isomorphism $i_{1}$, where

\begin{center}
$i_{1}:\mathcal{M}^{\ast }\rightarrow (S,\in )^{\mathcal{N}}.$
\end{center}

\noindent Finally, we verify that $S=N$. We first note that $S$ must be a
proper class in the sense of $\mathcal{N}$, since otherwise $\mathcal{N}$
would be able to define the satisfaction predicate for $(S,\in )^{\mathcal{N}%
}$, which coupled with the assumption that $\mathcal{N}$ is interpretable in 
$\mathcal{M}$, and $i:=i_{1}\circ i_{0}$ is an $\mathcal{M}$-definable
isomorphism between $\mathcal{M}$ and $(S,\in )^{\mathcal{N}}$, would result
in $\mathcal{M}$ being able to define a satisfaction predicate for itself,
which contradicts (an appropriate version of) \textit{Tarski's
Undefinability of Truth Theorem}\footnote{%
For a structure $\mathcal{M}$ let:
\par
$\mathrm{Th}^{+}(\mathcal{M})=\left\{ \left( \ulcorner \sigma \urcorner
,a\right) :\mathcal{M}\models \sigma (a)\right\} ,$ and $\ \mathrm{Th}^{-}(%
\mathcal{M})=\left\{ \left( \ulcorner \sigma \urcorner ,a\right) :\mathcal{M}%
\models \lnot \sigma (a)\right\} .$%
\par
\noindent With the above notation in mind, the version of Tarski's theorem
that is invoked here says that if\textit{\ }$\mathcal{M}$\textit{\ }is a
structure that interprets\textit{\ }$\mathsf{Q}$ and is endowed with a
pairing function, then\textit{\ }$\mathrm{Th}^{+}(\mathcal{M})$\textit{\ }and%
\textit{\ }$\mathrm{Th}^{-}(\mathcal{M})$\textit{\ }are\textit{\ }$\mathcal{M%
}$-inseparable, i.e., there is no $\mathcal{M}$-definable $D$\textit{\ }such
that $\mathrm{Th}^{+}(\mathcal{M})\subseteq D$ and $\mathrm{Th}^{-}(\mathcal{%
M})\cap D=\varnothing .$}.The transitivity of $S$ coupled with the fact that 
$S$ is a proper class in $\mathcal{N}$ together imply that $S$ contains all
of the ordinals of $\mathcal{N}$. Therefore, if $S\neq N,$ then arguing in $%
\mathcal{N}$, let $V_{\alpha }^{S}$ be $V_{\alpha }$ in the sense of $\left(
S,\in \right) $ and let

\begin{center}
$\alpha _{0}=$ the first ordinal $\alpha $ such that $V_{\alpha }=V_{\alpha
}^{S}$, but $V_{\alpha +1}\backslash V_{\alpha +1}^{S}\neq \varnothing $.
\end{center}

\noindent This makes it clear, in light of the assumption that $\mathcal{M}%
\trianglerighteq _{\mathrm{par}}\mathcal{N}$, and the fact that $i$ is an
isomorphism between $\mathcal{M}$ and $(S,\in )^{\mathcal{N}}$, that we have
a contradiction at hand since $\mathcal{M}$ believes that $\mathcal{N}$ sees
a `new subset' of $V_{i^{-1}(\alpha _{0})}$ of $\mathcal{M}$ that is missing
from $\mathcal{M}$. Hence $S=N$ and we may conclude that $i$ is an $\mathcal{%
M}$-definable isomorphism between $\mathcal{M}$ and $\mathcal{N}$.\hfill $%
\square $\medskip

\noindent \textbf{2.7.~Theorem. }$\mathsf{KM}$ \textit{is solid.}\medskip

\noindent \textbf{Proof.} Models of $\mathsf{KM}$ can be represented as
two-sorted structures of the form $\left( \mathcal{M},\mathcal{A}\right) $,
where $\mathcal{M}\models \mathsf{ZF}$; $\mathcal{A}$ is a collection of
subsets of $M$; and $\left( \mathcal{M},\mathcal{A}\right) $ satisfies the
full comprehension scheme. Suppose $\left( \mathcal{M},\mathcal{A}\right) $, 
$\left( \mathcal{M}^{\ast },\mathcal{A}^{\ast }\right) $, and $\left( 
\mathcal{N},\mathcal{B}\right) $ are models of $\mathsf{KM}$ such that:

\begin{center}
$\left( \mathcal{M},\mathcal{A}\right) \trianglerighteq _{\mathrm{par}%
}\left( \mathcal{N},\mathcal{B}\right) \trianglerighteq _{\mathrm{par}%
}\left( \mathcal{M}^{\ast },\mathcal{A}^{\ast }\right) ,$
\end{center}

\noindent and there is an $\left( \mathcal{M},\mathcal{A}\right) $-definable
isomorphism\textit{\ }

\begin{center}
$\widehat{i_{0}}:\left( \mathcal{M},\mathcal{A}\right) \rightarrow \left( 
\mathcal{M}^{\ast },\mathcal{A}^{\ast }\right) .$
\end{center}

\noindent As in the proof of Theorem 2.5 we note that $\widehat{i_{0}}$ is
naturally induced by $i_{0},$ where:

\begin{center}
$i_{0}:=\widehat{i_{0}}\upharpoonright M:\mathcal{M}\rightarrow \mathcal{M}%
^{\ast }$,
\end{center}

\noindent since $\widehat{i_{0}}(A)=\{i_{0}(m):m\in A\}$ for $A\in \mathcal{A%
}^{\ast }$.\medskip

\noindent $N$ forms a proper class in $\left( \mathcal{M},\mathcal{A}\right) 
$ since if $N$ forms a set, then so does $\mathcal{B}$, and $\widehat{i}_{0}$
is an $\left( \mathcal{M},\mathcal{A}\right) $-definable bijection between $%
M\cup \mathcal{A}$ and a subset of $N\cup \mathcal{B}$. Let $E:=\ \in ^{%
\mathcal{M}^{\ast }}$. Clearly $E$ is extensional. Furthermore, with the
help of $i_{0}$ and the assumption $\left( \mathcal{M},\mathcal{A}\right)
\trianglerighteq _{\mathrm{par}}\left( \mathcal{N},\mathcal{B}\right) $ it
is easy to see that $E$ is well-founded from the point of view of $\left( 
\mathcal{N},\mathcal{B}\right) $. The reader is asked to verify that an
argument very similar to the one used in the proof of Theorem 2.6 shows that 
$E$ is also set-like in the sense of $\left( \mathcal{N},\mathcal{B}\right) $%
. Theorem 2.6.2 can then be invoked to obtain an $\left( \mathcal{N},%
\mathcal{B}\right) $-definable isomorphism

\begin{center}
$k_{1}:$ $\mathcal{M}^{\ast }\rightarrow (S,\in ^{N})$
\end{center}

\noindent for some $\left( \mathcal{N},\mathcal{B}\right) $-definable
transitive $S\subseteq N$. The verification that $S=N$ is identical to the
corresponding part in the proof of Theorem 2.6 (and in particular uses
Tarski's undefinability of truth theorem). Let $k_{0}:=i_{0}^{-1}\circ
k_{1}^{-1}$. Clearly:\medskip

\noindent (5) $k_{0}:\mathcal{N}\rightarrow \mathcal{M}$ is an $\left( 
\mathcal{M},\mathcal{A}\right) $-definable isomorphism, and \medskip

\noindent (6) $k_{1}:\mathcal{M}^{\ast }\rightarrow \mathcal{N}$ is an $%
\left( \mathcal{N},\mathcal{B}\right) $-definable isomorphism. \medskip

\noindent Borrowing a notation from the proof of Theorem 2.5, let $\widehat{%
k_{0}}:\left( \mathcal{N},\mathcal{B}\right) \rightarrow \left( \mathcal{M},%
\mathcal{A}\right) $ be the natural extension of $k_{0}$, and $\widehat{k_{1}%
}:\left( \mathcal{M}^{\ast },\mathcal{A}^{\ast }\right) \rightarrow \left( 
\mathcal{N},\mathcal{B}\right) $ be the natural extension of $k_{1}$. Note
that both $\widehat{k_{0}}$ and $\widehat{k_{1}}$ are \textit{embeddings.}
Let $\widehat{k}:=\widehat{k_{0}}\circ \widehat{k_{1}}\circ \widehat{i_{0}}$%
; it is clear that:\medskip

\noindent (7) $\widehat{k}:\left( \mathcal{M},\mathcal{A}\right) \rightarrow
\left( \mathcal{M},\mathcal{A}\right) $ and $\widehat{k}$ is an $\left( 
\mathcal{M},\mathcal{A}\right) $-definable embedding.\medskip

\noindent Observe that (5) and (6), together with the assumption that $i_{0}$
is an isomorphism imply that the restriction $k:$ $\mathcal{M}\rightarrow 
\mathcal{M}$ of $\widehat{k}$ to `sets' is an automorphism of $\mathcal{M}$.
But since $\widehat{k}$ is $\left( \mathcal{M},\mathcal{A}\right) $%
-definable, $k(m)=m$ for all $m\in M$, thanks to the veracity of the scheme
of $\in $-induction\footnote{%
The scheme of $\in $-induction consists of the universal closures of
formulas of the form $\forall y\left( \forall x{\in }y\ \theta \left(
x\right) \rightarrow \theta (y)\right) \rightarrow \forall z\,\theta (z),$
where the parameters in $\theta $ are suppressed. It is easy to see that the
scheme of $\in $-induction is equivalent to the class-form of Foundation,
which asserts that every nonempty definable collection of sets has an $\in $%
-minimal element. The class-form of Foundation follows from the set-form of
Foundation and the comprehension scheme of $\mathsf{KM}$: suppose a class $C$
is nonempty, and let $\alpha _{0}$ be the first ordinal $\alpha $ such that $%
V_{\alpha }\cap C\neq \varnothing $. Then an $\in $-minimal member of $%
V_{\alpha _{0}}\cap C$ is also an $\in $-minimal member of $C$.} in $\mathsf{%
KM}$. This shows that $\widehat{k}$ is the identity map and in particular it
is surjective, which in turn implies that $\widehat{k_{0}}$ and $\widehat{%
k_{1}}$ are both surjective. This makes it clear that there is an $\left( 
\mathcal{M},\mathcal{A}\right) $-definable isomorphism between $\left( 
\mathcal{M},\mathcal{A}\right) $ and $\left( \mathcal{N},\mathcal{B}\right) $%
.\hfill $\square $ \medskip

Recall that $\widetilde{\mathsf{ZF}}$ was defined earlier, just before
Corollary 2.5.1.\medskip

\noindent \textbf{2.7.1.~Corollary.~}\textit{The following theory} $T$ 
\textit{is solid}:

\begin{center}
$T:=$ $\widetilde{\mathsf{ZF}}\backslash \{\mathsf{Power}\ \mathsf{Set}\}$ +
\textquotedblleft $\exists \kappa $ ($\kappa $ \textit{is strongly
inaccessible}, \textit{and }$\forall x\ \left\vert x\right\vert \leq \kappa $%
\textquotedblright $).$
\end{center}

\noindent \textbf{Proof.~}In light of Remark 2.1.1, this follows from
Theorem 2.8 and the well-known bi-interpretability of $T$ with $\mathsf{KM}%
+\Pi _{\infty }^{1}$-$\mathsf{AC}$, where $\Pi _{\infty }^{1}$-$\mathsf{AC}$
is the scheme of Choice.\footnote{%
This bi-interpretability was first noted by Mostowski; a modern account is
given in a recent paper of Antos \& Friedman \cite[section 2]{Antos-Friedman}%
, where $\mathsf{KM}+\Pi _{\infty }^{1}$-$\mathsf{AC}$ is referred to as $%
\mathsf{MK}^{\ast }$, and $T$ is referred to as $\mathsf{SetMK}^{\ast }$.}%
\hfill $\square $\medskip

\noindent \textbf{2.8.~Remark.~}An examination of the proofs in this section
make it clear that for each positive integer $n$, the theories $\mathsf{Z}%
_{n}$ ($n$-th order arithmetic) and $\mathsf{KM}_{n}$ ($n$-th order
Kelley-Morse theory of classes) are solid theories (where $\mathsf{Z}_{1}:=%
\mathsf{PA}$, and $\mathsf{KM}_{1}:=\mathsf{ZF}$). This observation, in
turn, implies that the theory of types $\mathsf{Z}_{\omega \text{ }}$ (with
full comprehension) whose level-zero objects form a model of $\mathsf{PA}$
(equivalently $\mathsf{ZF}_{\mathsf{fin}}+\mathsf{TC)}$, and the theory of
types $\mathsf{KM}_{\omega }$ whose level-zero objects form a model of $%
\mathsf{ZF}$ are also solid theories. Thus, the list of theories whose
solidity is established in this section can be described (up to
bi-interpretability) as $\left\{ \mathsf{Z}_{n}:1\leq n\leq \omega \right\}
\cup \left\{ \mathsf{KM}_{n}:1\leq n\leq \omega \right\} .$\bigskip 

\section{Examples and Questions}

All of the theories $T$ whose solidity was established in Section 2 are
sequential\footnote{%
A sequential theory is a theory that has access to a definable `$\beta $%
-function'\ for coding finite sequences of objects in the domain of
discourse.} theories which have an interpretation $\mathbb{N}$ for `numbers'
for which the full scheme $\mathsf{Ind}_{\mathbb{N}}(\mathcal{L}_{T})$ of
induction is $T$-provable, so one may ask whether the $T$-provability of $%
\mathsf{Ind}_{\mathbb{N}}(\mathcal{L}_{T})$ within a sequential theory is a
sufficient condition for solidity. A simple counterexample gives a negative
answer: let $\mathsf{PA}(\mathsf{G})$ be the natural extension of $\mathsf{PA%
}$ in which the induction scheme is extended to formulae in the language
obtained by adding a unary predicate $\mathsf{G}$ to the language of
arithmetic. To see that $\mathsf{PA}(\mathsf{G})$ is not solid, consider the
extensions $T_{1}$ and $T_{2}$ of $\mathsf{PA}(\mathsf{G})$, where:

\begin{center}
$T_{1}:=\mathsf{PA(G)}+\forall x(\mathsf{G}(x)\leftrightarrow x=1)$ and $%
T_{2}:=\mathsf{PA(G)}+\forall x(\mathsf{G}(x)\leftrightarrow x=2).$
\end{center}

\noindent Clearly the deductive closures of $T_{1}$ and $T_{2}$ are
distinct, and yet it is easy to see that $T_{1}$ and $T_{2}$ are
bi-interpretable. This shows that $\mathsf{PA(G)}$ is not tight, and
therefore not solid.\medskip

With the help of \cite[Theorem 4.9 \& Remark 4.10]{ESV} one can also show
that the theory $\mathsf{ZF}_{\mathsf{fin}}$ is not tight, even though as
shown in Corollary 2.2.1 its strengthening by \textsf{TC} is a solid theory.
Another example of a theory that fails to be tight is $\mathsf{ZF\backslash
\{Foundation\}}$. To see this, consider $T_{1}:=\mathsf{ZF}$, and

\begin{center}
$T_{2}:=\mathsf{ZF\backslash \{Foundation\}}+\exists !x(x=\{x\}+\forall
t\exists \alpha (t\in V_{\alpha }(x)),$
\end{center}

\noindent where $\alpha $ ranges over ordinals, and

\begin{center}
$V_{0}(x):=x$, $V_{\alpha +1}(x):=\mathcal{P}(V_{\alpha }(x)),$ and $%
V_{\alpha }(x):=\bigcup\limits_{\beta <\alpha }V_{\beta }(x)$ for limit $%
\alpha .$
\end{center}

\noindent Then $T_{1}$ and $T_{2}$ are extensions of $\mathsf{ZF\backslash
\{Foundation\}}$ with distinct deductive closures, and yet, the
bi-interpretability of $T_{1}$ and $T_{2}$ can be established by well-known
methods: the relevant interpretations are $\mathcal{I}$ and $\mathcal{J}$,
where $\mathcal{I}$ is the classic von Neumann interpretation of $\mathsf{ZF}
$ in $\mathsf{ZF\backslash \{Foundation\}}$, and $\mathcal{J}$ is the
classic Rieger-Bernays interpretation that adds a single `Quine atom' (i.e.,
a set $s$ such that $s=\{s\}$) to a model of $\mathsf{ZF}$.\medskip

However, we do not know whether $T\vdash \mathsf{Q}^{\mathbb{N}}+\mathsf{Ind}%
_{\mathbb{N}}(\mathcal{L}_{T})$ for every solid sequential theory (for an
appropriate choice of numbers $\mathbb{N}$). This motivates the following
question, since by a general result of Montague \cite{Montague} the $T$%
-provability of $\mathsf{Q}^{\mathbb{N}}+\mathsf{Ind}_{\mathbb{N}}(\mathcal{L%
}_{T})$ implies that $T$\ is not finitely axiomatizable.\medskip

\noindent \textbf{3.1.~Question.~}\textit{Is there a consistent sequential
finitely axiomatized theory that is solid?}\medskip

\noindent The question below arises from reflecting on the results of
Section 2 and noting that the proofs of solidity of each of the theories $T$
established in Section 2 uses the `full power'\ of $T$.\medskip

\noindent \textbf{3.2.~Question.~}\textit{Is there an example} $T$\textit{\
of one of the theories whose solidity is established in Theorem 1.2}, 
\textit{and some \textbf{solid} }$T_{0}\subseteq T$ \textit{such that the
deductive closure of} $T_{0}$ \textit{is a proper subset of the deductive
closure of }$T?$

\section{Acknowledgements}

It is a pleasure and an honor to present this paper in a volume that
celebrates Albert Visser's scholarship; I am grateful to Albert for bringing
his Theorem 1.1 to my attention. Thanks also to Andr\'{e}s Caicedo and Radek
Honz\'{\i}k, whose interest in the $\mathsf{ZF}$-case of Theorem 1.2
provided additional impetus for writing up the results here; and to the
anonymous reviewer for invaluable help in weeding out infelicities of an
earlier draft. Hats off to Jan, Joost and Rosalie for their dedication in
bringing this volume to fruition.

\end{document}